\documentclass [11pt] {amsart}
\usepackage{amsthm,amsmath,amssymb,amscd,graphics,enumerate}
\usepackage{latexsym}
\usepackage{verbatim}
\usepackage{enumerate}
\usepackage{amsthm}
\usepackage{amsmath}
\usepackage{amscd}

\newcounter{theorem}[section]
\numberwithin{equation}{section}

\newtheorem{Thm}[theorem]{Theorem}
{\theoremstyle{remark}
 
\newtheorem{Rem}[theorem]{\text{\textbf{Remark}}} }
\newtheorem{Def}[theorem]{Definition}
\newtheorem{Lem}[theorem]{Lemma}
\newtheorem{Prop}[theorem]{Proposition}
\newtheorem{Cor}[theorem]{Corollary}

\newtheorem{assumptions}[theorem]{Assumptions}

{\theoremstyle{definition}
 }

\newcommand{\RR}{{\mathbb{R}}}

\newcommand{\CC}{{\mathbb{C}}}
\newcommand{\QQ}{{\mathbb{Q}}}

\newcommand{\PP}{{\mathbb{P}}}
\newcommand{\ZZ}{{\mathbb{Z}}}

\newcommand{\mO}{\mathcal{O}}
\newcommand{\mX}{\mathcal{X}}

\newcommand{\mC}{\mathcal{C}}

\newcommand{\mE}{\mathcal{E}}

\newcommand{\mK}{\mathcal{K}}

\title{Note on characterization of projective spaces}
\author{Jiun-Cheng Chen}
\address{Department of Mathematics\\ Northwestern
University\\ 2033 Sheridan Road\\ Evanston\\ IL 60208-2370\\ USA}
\email{jcchen@math.northwestern.edu}

\author{Hsian-Hua Tseng}
\address{Department of Mathematics\\ University of British Columbia\\ 1984 Mathematics Road\\ Vancouver\\ B.C. V6T 1Z2\\ Canada}
\email{hhtseng@math.ubc.ca}

\date{\today}

\begin{document}
\begin{abstract}
We prove a numerical characterization of $\PP^n$ for varieties with at worst  isolated local complete intersection quotient singularities. In dimension three, we prove such a numerical characterization of $\PP^3$ for normal $\QQ$-Gorenstein  projective varieties. 
\end{abstract}
\maketitle
\section{Introduction}
We work over $\CC$. Recall the following numerical characterization of the
projective space $\PP^n$.
\begin{Thm}[\cite{cmsb02} and \cite{ke01}]
Let $X$  be a smooth  projective variety of dimension $n \geq 3$. Assume that
$C \cdot (-K_X) \geq n+1$ for all curves  $C \subset X$.
Then $X\cong \PP^n$.
\end{Thm}
This result  was first proved by Cho, Miyaoka and Shepherd-Barron
\cite{cmsb02} and later simplified  by Kebekus  \cite{ke01}. The main
goal of this paper is to relax the assumption on smoothness.  

\begin{Thm}\label{main} Let $X$ be a projective variety
of dimension $n \geq 3$ with at most isolated local complete intersection quotient (LCIQ) singularities.
Assume that there is a $K_X$-negative extremal ray $R$ such that $C\cdot (-K_X) \geq n+1$ for every curve $[C] \in R$. 
Then $X \cong\PP^n$.
\end{Thm}
Note that the numerical condition in Theorem \ref{main} is weaker: we only require this condition only for curves in one extremal ray, instead of all curves.

The next corollary follows immediately.
\begin{Cor}\label{main2}
Let $X$  be a projective variety of dimension $n \geq 3$ with at
most isolated LCIQ singularities. If $C \cdot (-K_X) \geq n+1$ for
all curves $C \subset X$, then $X \cong \PP^n$.
\end{Cor}
Combining with methods from the minimal model program (MMP), we obtain the following 
stronger result when $dim\; X=3$.
\begin{Thm}\label{threefold}
Let $X$ be a normal $\QQ$-Gorenstein  projective variety of dimension $3$ such  that $C \cdot (-K_X) \geq 4$ for every curve
$C \subset X$. 
Then $X \cong \PP^3$.
\end{Thm}

 Our  main strategy is to produce enough  (minimal) rational curves in $X$.  For smooth $X$, deformation theoretic methods can be used to achieve this goal. For singular $X$, it is difficult to produce rational curves  via  deformation theory of maps from $\PP^1$ to $X$.  However, if $X$ has only local complete intersection quotient (LCIQ)  singularities, we may proceed by considering an alternative:  Instead of maps from $\PP^1$ to $X$, we consider representable maps from a {\em twisted curve}\footnote{Roughly speaking, this is a one dimensional Deligne-Mumford stack with isolated cyclic quotient stack structures and with coarse moduli space being a nodal curve \cite{av02}.} to the stack $\mX$ whose coarse moduli space is $X$. We obtain a lower bound on the dimension of  the space of twisted stable maps expressed in terms of the $-K_X$-degree and the number of twisted points, see \cite{ct05}.
 

In dimension $3$, we can assume much less on singularities. 
When the threefold  $X$ has possibly worse singularities, the deformation theoretic methods do not apply directly.
A  natural strategy is to find a suitable partial resolution $Z$ of $X$ and try to prove $Z \simeq \PP^3$
via deformation theory of twisted stable maps. Since $Z \simeq \PP^3$, the morphism $h$ has to be an isomorphism  (and hence  $X \simeq \PP^3 $).   
However, the numerical condition $C \cdot (-K_X) \geq 4\; \forall \; C \subset X$ is not stable under  birational  modification;   
this is the main motivation of the weaker condition in Theorem \ref{main}.

We now sketch our proof of Theorem \ref{threefold}. Take a terminal $\QQ$-factorial model  $Z$ of $X$. We can prove that the stronger numerical condition 
$$ C \cdot (-K_X) \geq 4 \; \forall \;C \subset X$$
implies  a weaker numerical condition (needed in Theorem \ref{main}): there is a $K_Z$-negative  extremal ray  $R$ in $\overline{NE}(Z)$ such that  $C \cdot (-K_Z) \geq 4\; \forall  \;C \subset Z \text{ and } [C] \in R$.  We then prove 
 Theorem \ref{threefold} by applying  Theorem \ref{main} to $Z$.
  
The rest of this paper is organized as follows. In Section \ref{vfree} we collect  basic facts on free and very free curves. When $X$ is smooth, all corresponding
statements in that section can be found in \cite{ko96}. The proofs
for  the singular cases require few (if any) changes since we make a
very strong assumption that a general rational curve from a suitable
family does not meet the singular locus $X_{sing}$ of $X$. This
assumption holds for varieties we are interested in. The proof  is
in Section~\ref{untwistedsec}. The main novelty is to apply
techniques developed in \cite{ct05}. In Section~\ref{kebekus}, we
recall results from \cite{ke00}, \cite{ke01} and prove Theorem \ref{main}. We do not claim any
originality of  results in that section; they are included in this
paper for  completeness of the proof and for the reader's
convenience. In the last section, we present a proof of Theorem \ref{threefold}.
\section*{Acknowledgments} We thank Stefan Kebekus for
helpful discussions. Part of this research was conducted  while we were
attending  the AMS summer institute of algebraic geometry. We like to
thank the conference organizers. The first author started working on
this paper during his visit at National Center of Theoretical
Sciences in Hsinchu, Taiwan. He likes to thank their hospitality and
support. The research of the second author is
supported in part by Clay Mathematics
Institute Liftoff fellowship program.
\section{Free and very free  curves}\label{vfree}
We recall several basic facts on  free and very free rational curves
in this section. When the variety $X$ is smooth, all statements we
collect here can be found in \cite{ko96}. These statements still
hold for singular $X$ as long as the rational curves do not meet\footnote{This condition  is  quite strong. It is easy to find examples that
all rational curves have to meet the singular locus $X_{sing}$.} the
singular locus $X_{sing}$.

 Let $f: \PP^1 \to X$
be a rational curve. Assume that $f: \PP^1 \to X$ is birational to
its image and $f(\PP^1)$ does not meet the singular locus
$X_{sing}$. Since $X$ is smooth along $f(\PP^1)$, the pull-back of
the tangent bundle $TX$ to $\PP^1$ is a vector bundle. It is  a
well-known
 result that $$f^*TX \cong \mO(a_1) \oplus \mO(a_2) \oplus \cdots \oplus
\mO(a_n)$$ (this is usually attributed to Grothendick, however, it
seems has been known before Grothendick) where $a_1 \geq a_2 \geq
\cdots \geq a_n$ are integers.

\begin{Def} Notation and assumptions as
above. The rational curve $f:\PP^1 \to X$  is free if $f^*TX$ is
generated by global sections. The rational curve is very free if
$f^*{TX} \otimes I_{p}$ is generated by global section  $\forall\; p
\in \PP^1$.
\end{Def}
\begin{Rem}
It is clear that a rational curve is  free (res. very free) if and only if  $a_i  \geq 0$ (res. $a_i \geq 1$).
\end{Rem}
The next proposition is a standard result. It only  works when the ground filed
$k=\bar{k}$ is uncountable and $ch(k)=0$, e.g.
$k=\CC$.
\begin{Prop}[\cite{ko96} Theorem II.3.11]\label{free}
Let $f:\PP^1 \to X$ be any rational curve which does not meet
$X_{sing}$. There exists a subset $X^{free} \subset X$,  which is  the
complement of  a countable union of proper subvarieties of $X$, such
that if $f(\PP^1)$  intersects $X^{free}$, then $f:\PP^1 \to X$ is
free.
\end{Prop}
We usually need to work on a family of (non-singular) rational curves. Let $Z$ be
an irreducible variety  and
$$
\begin{CD}
U @ >{\operatorname{F}}>>X\\
@V{\operatorname{\pi_{Z}}}VV \\
Z
\end{CD}
$$
 a family of morphisms over $Z$ whose domain curves are
non-singular rational curves. Let $[f] \in Z$ be a member of this
family.

\begin{Lem}[\cite{ko96} Theorem II.3.14.3]\label{embedding}
Notation as above. If $a_3 \geq 1$, then a general deformation $g:
\PP^1 \to X$ of $f$ is an embedding.
\end{Lem}
\begin{Prop}[\cite{ko96} Corollary II.3.10.1]\label{ample}
Notation as above. Let $x$ be a general point of $X$. Assume that
$g(0)=x$ for all $[g] \in Z$. Let $[f] \in Z$ such that  the image
of $f$ does not meet $X_{sing}$. Assume further that $F: U \to X$ is
dominant. Then for a general member $[h] \in Z$, the image
$h(\PP^1)$ does not meet $X_{sing}$ and  the vector bundle $h^* TX$
is ample.
\end{Prop}

\section{Untwisted rational curves}\label{untwistedsec}
We assume that $X$ is a projective variety of dimension $n \geq 3$
with at most isolated LCIQ singularities in this section. For such
an $X$, one can find a proper LCI Deligne-Mumford stack $\pi: \mX
\to X$ such that $X$ is a coarse moduli space of $\mX$ and  $\pi$ is
an isomorphism over $X_{reg}=X \setminus X_{sing}$. It is not hard
to see that $K_{\mX}=\pi^{*}K_X$ and $$\mC \cdot K_{\mX}= C
\cdot K_X$$ for any (twisted) curve $\mC \to \mX$ with coarse curve $C$ \cite{ct05}.

The main goal of this section is to prove the next proposition.

\begin{Prop}\label{general2}
Notation and assumptions as above. Assume further that there
is a $K_X$-negative extremal   ray $R \subset \overline{NE}(X)$ such
that $C \cdot (-K_X) \geq n+1$ for
 every curve $[C] \in R$.
Let $x \in \mX$  be a general point. Then there is a rational curve
$f:\PP^1 \to \mX$ which passes through $x$ and does
not meet $\pi^{-1}(X_{sing})$.
\end{Prop}

We start with the following  proposition.
\begin{Prop}\label{lessthan}
Notation as above.  Let $R$ be any $K_X$-negative extremal ray. Then
there is a twisted stable map $f:\mC \to \mX$ such that $\mC$ has at
most one stacky point, $-K_X \cdot \mC \leq n+1$ and  $[\pi \circ f
(\mC)] \in R$.
\end{Prop}
\begin{proof}
Let $R$ be such an extremal ray and $f:C \to X$  a curve whose
class  is in the ray $R$. Choose a finite cover $D \to C$ such that
$D \to X$ can be lifted to $g:D \to \mX$. We may assume that $\pi:
\mX \to X$ and $g:D \to \mX$ are  defined over $\ZZ \subset S \subset
\CC$ where $S$ is finitely generated  over $\ZZ$.
The mod $p$ reductions of $\pi:\mX \to X$ and $g:D \to \mX$
are denoted by  $\pi_p:\mX_p \to X_p$ and $g_p:D_p \to \mX_p$ respectively.
It follows from \cite{ct05} that $R$ remains extremal for almost all
$p$. The  characteristic $p$ case of this proposition  follows  from \cite{ct05}: essentially,
this follows by applying the Frobenius map and  bend-and-break techniques.
Let $h_p: \mC_p  \to \mX_p$ be the twisted curve satisfying
the conditions
\begin{enumerate}
\item $\mC_p$ has at most $1$ twisted point,
\item $-K_{X_{p}} \cdot \mC_p \leq n+1$, and
\item $[\pi_p \circ h_p(\mC_p)] \in R$.
\end{enumerate}
as  in the proposition.  Also note that the coarse moduli space $C_p$ of $\mC_p$ is rational.

To prove the  characteristic $0$ case, we apply the the standard technique
of reduction mod $p$ (\cite{ko96} Meta Theorem II.5.10.3).

Consider the two conditions  for a twisted stable map $h:\mC \to
\mX$ (over the base scheme $Spec\;S$)
\begin{enumerate}
\item $-K_X \cdot \mC \leq n+1$, and
\item $[\pi \circ h(\mC)] \in R$.
\end{enumerate}
There are only  finitely many components  of twisted stable maps
such
 that the corresponding twisted stable maps
satisfy these two conditions. Note that the locus where the source
curve has at most $1$ stacky point is open. The image of these
components in $Spec\;S$ is constructible.
  Since the image in $Spec \;S$
contains almost all $p$, it also contains the generic fiber. This
shows that there is a  twisted stable map with at most $1$ stacky
point.
\end{proof}

Unless mentioned  otherwise,
we make further assumptions in the rest of this paper.
\begin{assumptions}\label{assu} From now on we assume that $X$ is a
projective variety of $dim X=n \geq 3$ with isolated LCIQ
singularities and has the property that there is a $K_X$-negative
 extremal ray $R$ such that
 $C \cdot (-K_X) \geq n+1$ for
 every curve $[C] \in R$.
\end{assumptions}

Let $R$ be an extremal ray as in  Assumptions~\ref{assu}. By
Lemma~\ref{lessthan}, there is a twisted stable map $f: \mC \to \mX$
such that $\mC \cdot (-K_{\mX}) \leq n+1$. Together with
Assumptions~\ref{assu}, it follows that $\mC \cdot (-K_{\mX})=n+1$.

Recall that a twisted stable map $f:\mC \to \mX$ is an equivalence
class of twisted curves; that is, we forget the parametrization of
$\mC$. Choosing a representative of this equivalence class yields a
twisted curve into $\mX$. (Usually it does not matter which
representative we choose.) Abusing the notation, we still  use $f:
\mC \to \mX$ to denote this twisted curve.

The next lemma is a simple observation.
\begin{Lem}\label{bir}
Assumptions as above. Let $f: \mC \to \mX$ be a twisted rational
curve with at most $1$ twisted point on $\mC$. Assume further that
the class of the coarse curve $[\overline{\pi \circ f}(\PP^1)] \in
R$ and $\mC \cdot (-K_{\mX}) =n+1$.
 Then $f: \mC \to \mX$ is birational onto its image.
\end{Lem}
\begin{proof}
This follows from the fact the $\mC \cdot (-K_\mX)=n+1$ is the
smallest possible intersection number by Assumptions~\ref{assu}.
Thus $f:\mC \to \mX$ can not be a multiple cover of its image.
\end{proof}

Let $f:\mC \to \mX$ be a twisted stable map with at most $1$ twisted
point and $\mC \cdot (-K_\mX)=n+1$. Let $y \in f(\mC) \subset \mX$
be a general point in the image of $\mC$ (we do not and can not
assume $y$ is general in $\mX$ for the time being).
Here we abuse the notation again: we use $f:\mC \to \mX$
to denote a representative in the equivalence class. We may  assume
$f(0)=y$.
\begin{Prop}\label{dim} Let $\mC  \to \mX$ be a twisted stable map with at most $1$ twisted
point and $\mC \cdot (-K_\mX)=n+1$. Then $dim_{[f]} Mor(\mC,
\mX,n+1, f\mid_{0}) \geq 2$.
\end{Prop}
\begin{proof}
We start with the case when $\mC$ does have a twisted point, denote
it by $\infty$. Consider the dimension estimate from \cite{ct05}:
 $$dim_{[f]} Mor(\mC, \mX, f \mid_{0}) \geq \mC \cdot (-K_{\mX})+
  (1-1)(dim\; \mX)-(\text{age term at $\infty$}).$$
  Note that the age term  is (strictly) less than $n$
  (see  the proof of Lemma4.6 in \cite{ct05}).
It follows that $dim_{[f]} Mor(\mC, \mX,n+1, f \mid_{0})> n+1-n=1$.
Since $dim_{[f]} Mor(\mC, \mX,n+1, f \mid_{0})$ is an integer, we
get $$dim_{[f]} Mor(\mC, \mX,n+1, f \mid_{0}) \geq 2.$$

The case when $\mC$ has no twisted point, i.e. $\mC \cong \PP^1$,
follows from a similar dimension estimate (without  the age term)
$$dim_{[f]} Mor(\PP^1, \mX, f \mid_{0}) \geq \PP^1 \cdot (-K_{\mX})+
  (1-1)(dim\; \mX)=n+1 \geq 2.$$
\end{proof}
\begin{Rem}\label{inf}
\hfill
\begin{enumerate}
\item
Suppose that $\mC$ does have a twisted point, denoted by $\infty$.
Note that any (representable) automorphism $\sigma:\mC \to \mC$ will
fix the twisted point $\infty \in \mC$. It is not hard to see that
 $dim\; Aut\,(\mC,0) =1$ where $Aut\,(\mC,0)$ consists of automorphisms
which fix $0$. Intuitively,  an automorphism $\sigma \in
Aut\,(\mC,0)$ can be viewed as an automorphism of $\PP^1$ fixing $0$
and $\infty$.

\item One  observes that $g(\infty)=f(\infty)$ when $[f]$ and $[g]$
are in the same component of $Mor(\mC, \mX, f \mid_{0})$ since the
image of a twisted point has to be a twisted point on $\mX$, and
these points are isolated by assumption.
\end{enumerate}
\end{Rem}

\begin{Lem}\label{untwisted}
Notation as above. The curve $\mC$ has no twisted point, i.e. $\mC
\cong \PP^1$.
\end{Lem}
\begin{proof}
Note that $\mC \cdot (-K_{\mX})=n+1$ is the smallest possible
intersection number for curves whose classes are in $R$. Therefore
the curve $f: \mC \to \mX$ can not be further broken. Suppose  that
$\mC$ does have a twisted point. By Remark~\ref{inf}, the image of
the twisted point $\infty$ is fixed for every $g:\mC \to \mX$ in the
same component. The dimension of such morphisms is at most $1$, i.e.
the dimension of  $Aut(\mC,0)$; otherwise we can apply bend and
break to get a curve with smaller $-K_{\mX}$ degree. However,
$dim_{[f]} Mor(\mC, \mX,n+1, f \mid_{0}) \geq 2$ by
Proposition~\ref{dim}. This shows that  the curve $\mC$ has no
twisted point at all.
\end{proof}

\begin{Cor}\label{dim1}
Notation as above. $dim_{[f]} Mor(\PP^1, \mX,n+1, f \mid_{0}) \geq
n+1$.
\end{Cor}
\begin{Rem}\label{dim2}
\hfill
\begin{enumerate}
\item The equality in Corollary~\ref{dim1} holds whenever
$f:\PP^1 \to \mX$ is very free.
\item
It's easy to compute the dimension of the space of twisted stable
maps passing through a point $x\in \mX$. Abusing the notation, we consider $f:\PP^1
\to \mX$ as a twisted stable map. Let $[f] \in H_x \subset
\mK(\PP^1,\mX,n+1)$ be an irreducible component of the subfamily of twisted
stable maps whose images contain $x$. It follows easily that
$$dim\;_{[f]} H_x=dim_{[f]} Mor(\PP^1, \mX,n+1, f
\mid_{0})-dim\;Aut(\PP^1,0)$$ $$ \geq n+1-2=n-1.$$
\end{enumerate}
\end{Rem}

Let $[f] \in M \subset  Mor(\PP^1, \mX,n+1, f \mid_{0})$ be an
irreducible component of dimension at least $n+1$. (In fact, there
is only one irreducible component which contains $[f]$.)
\begin{Cor}\label{general}
Notation as above. A general member in $M \subset Mor(\PP^1, \mX,n+1, f
\mid_{0})$ does not meet the preimage of the singular locus of $X$.
\end{Cor}
\begin{proof}

Consider the  subfamily ($\subset M$) of morphisms  which  meet
 $\pi^{-1}(X_{sing})$.
We claim this family  has dimension at most $2$:
Let $t \in \PP^1$ and $z \in \pi^{-1}(X_{sing})$.
The dimension of the  subfamily
of morphisms $h:\PP^1 \to\mX$
 such that $h(0)=y$ and $h(t)=z$ is at most $1$ by bend and break.
Since $t$ is chosen from a $1$-dimensional
family, i.e. $\PP^1$,
and $\pi^{-1}(X_{sing})$ is finite, the claim follows.
The corollary follows easily from our claim.
\end{proof}

\begin{proof}[Proof of Proposition~\ref{general2}]
Choose a general
$[f] \in Mor(\PP^1, \mX,n+1, f\mid_{0})$
such that $f(\PP^1)$ does not
meet $\pi^{-1}(X_{sing})$ (Corollary~\ref{general}).
We may also assume that
$$Mor(\PP^1, \mX,n+1, f\mid_{0})$$ is smooth
at $[f]$ by Proposition~\ref{ample} (\cite{ko96} Corollary
II.3.10.1). Take the irreducible component $M \subset  Mor(\PP^1,
\mX,n+1, f\mid_{0})$ which contains $[f]$ (there is only one such
component). Let
$$\begin{CD}
U_{M}@ >{F_{M}}>>X\\
@V{\pi_{M}}VV \\
M
\end{CD}$$
be the family of morphisms into $X$ over $M$. Note that the image of
a general member $[f] \in M$ does not meet $\pi^{-1}(X_{sing})$ (see
Corollary~\ref{general}). To conclude the proof, it suffices to show
that $F_M$ is dominant. We prove this by a simple dimension count.
Let $x_1 \in X$ be a general point. Consider the fiber of $F_ M$
over $x_1$. The dimension of $F_M^{-1}(x_1)$ is at most $2$;
otherwise we can apply bend and break to obtain a lower
$(-K_X)$-degree rational curve in the ray $R$. Note that $dim_{[f]} M=
n+1$ (Corollary~\ref{dim1} and Remark~\ref{dim2}) and
$dim\;U_M=n+2$. This shows that $F_M$ is dominant.
\end{proof}

It is well-known \cite{av02} that  $\mK(\PP^1,\mX,n+1)$, the space (stack) of
twisted stable maps with $(-K_X)$-degree $n+1$, is quasi-finite over
 the coarse moduli space
$K(\PP^1,X,n+1)$.
 Since $n+1$ is the minimal $(-K_X)$-degree
 by assumption, the domain curves do not degenerate and
the universal family of stable map $U \to   \mK(\PP^1, \mX,n+1)$ is
a $\PP^1$-bundle over $\mK(\PP^1, \mX,n+1)$.
\begin{Lem}
The  proper Deligne-Mumford stack $\mK(\PP^1, \mX,n+1)$ is a
projective scheme.
\end{Lem}
\begin{proof}
By Lemma \ref{bir}, any stable map $[f] \in \mK(\PP^1, \mX,n+1)$ is
birational to its image in $\mX$. Therefore, $f: \PP^1 \to \mX$ has
no nontrivial automorphism. This shows  $\mK(\PP^1, \mX,n+1)$ is a
proper algebraic space. Since it is also quasi-finite over the
projective scheme $K(\PP^1,X,n+1)$, it is  a projective scheme.
\end{proof}

Let $[f] \in \mK(\PP^1, \mX,n+1)$  be a twisted stable map into
$\mX$ which does not meet the preimage of $X_{sing}$. Let $Z \subset
\mK(\PP^1,X ,n+1)$ be an irreducible component  which contains $[f]$
and $\tilde{Z}$ the normalization. Consider the finite  morphism
$\mK(\PP^1,\mX,n+1) \to K(\PP^1,X,n+1)$. Let $Z'$ be the component
of $K(\PP^1,X,n+1)$ which contains the  image of $Z$, i.e. stable
maps  of the form $[\pi \circ h]$ with $[h] \in Z$. We also take the
normalization $\tilde{Z}$ of $Z$.
The next lemma compares these two components. It is needed in the next
section.
\begin{Lem}\label{iso}
The natural map $\tilde{Z} \to \tilde{Z'}$ is an
isomorphism.
\end{Lem}
\begin{Rem}
We do not claim that
$\mK(\PP^1,\mX,n+1) \cong K(\PP^1, X, n+1)$ at this moment, though it turns out
to be the case.
We have not excluded the possibility that
 there is a  component  of $K(\PP^1, X, n+1)$ consisting of rational curves
 which can not be lifted to $\mX$.
\end{Rem}
\begin{proof}
The morphism  $\tilde{Z} \to \tilde{Z'}$ is finite. Since both
$\tilde{Z}$ and $\tilde{Z'}$ are normal, it suffices to show the
morphism is also birational. This is clear since a general twisted
stable map  $[h] \in \tilde{Z}$ does not meet the set
$\pi^{-1}(X_{sing})$. That is, $h:\PP^1 \to \mX$ lies in
$\pi^{-1}(X_{reg})$. Since $\pi: \mX \to X$ is an isomorphism over
$X_{reg}$,  we can identify an open set of $\tilde{Z}$ with an open
set of $\tilde{Z'}$. This concludes the proof.
\end{proof}

\section{Results from \cite{ke00} and \cite{ke01}}\label{kebekus}
In this section, we return to work on the variety $X$, rather than
the stack $\mX$. All results in this section are taken from
\cite{ke00} and \cite{ke01}. We follow his notation as closely as
possible. Let $x$ be a general point on $X$. Let $f:\PP^1 \to
\mX$ be a twisted rational curve such that
\begin{enumerate}
\item $\PP^1 \cdot _{f} (-K_{\mX})=n+1$,
\item $[\pi \circ f (\PP^1)] \in R$, and
\item the image $f(\PP^1)$ contains $x$ and does not
meet $\pi^{-1}(X_{sing})$.
\end{enumerate}
 The existence of such a curve follows
from Corollary~\ref{general2}. Denote by $[f] \in H_x \subset
\mK(\PP^1, \mX,n+1)$ an irreducible component of the subfamily of
stable map through $x$. Recall that $dim_{[h]} H_x=n-1$ for a
general $[h] \in H_x$ (Remark~\ref{dim2}). Let $\tilde{H}_x \to H_x$
be the normalization. By Lemma~\ref{iso}, the variety $\tilde{H}_x$
can be viewed as an irreducible component of the subfamily of stable
maps (into $X$) passing through $x$. Consider the diagram
$$
\begin{CD}
U_x @>{i_x}>>X\\
@V{\pi_x}VV\\
\tilde{H}_x
\end{CD}$$
where $\pi_x: U_x \to \tilde{H}_x $ is  the universal family over
$\tilde{H}_x$ and $i_x: U_x \to X$ the universal
morphism into $X$.
Let
 $\tilde{X}  \to X$ be the blow up of
 $X$ at the general point $x$. There is a rational
map $\tilde{i}_x: U_x \dashrightarrow  \tilde{X}$ lifting the
morphism $i_x: U_x \to X$.  We need to consider $\tilde{H}_x$ since
$H_x$ may not be normal a priori. Taking normalization  is not
necessary in \cite{ke01} since he proves $H_x$ is smooth for a very
general $x$ at the very beginning.

 We list some properties of $\tilde{H}_x$.
\begin{Prop}
\hfill
\begin{enumerate}
\item The variety $\tilde{H}_x$ is projective and  is
smooth away from a finite set of points.
\item The evaluation
morphism $i_x$ is finite away from $i^{-1}(x)$. In particular, the
morphism $i_x$ is surjective.
\item For a general
point $[l] \in \tilde{H}_x$ the corresponding curve $l \subset X$
does not meet the preimage of $X_{sing}$.
\item For a general point $[l] \in \tilde{H}_x$ the
corresponding curve $l \subset X$ is smooth and the restriction of
the tangent bundle $TX \mid_{l}$  is ample.
\end{enumerate}
\end{Prop}
\begin{proof}
The variety $H_x$ is non-empty by Proposition~\ref{general2}.  The
first half of (1) follows easily from \cite{av02} and the fact that
these curves do not degenerate. The smoothness statement follows
from Proposition~\ref{ample} (\cite{ko96} Corollary II.3.10.1). (2)
follows from the standard bend-and-break techniques, \cite{ko96}
Corollary II.5.5. (3) is  just a restatement of
Corollary~\ref{general}. The first part of (4) follows from
Lemma~\ref{embedding} (\cite{ko96} Theorem II.3.14); the second part
follows from Proposition~\ref{ample} (\cite{ko96} Corollary
II.3.10.1).
\end{proof}
The next proposition can be found in \cite{ke00} or \cite{ke01}.
\begin{Prop}[=\cite{ke01} Corollary 2.3]
The preimage $i^{-1}(x)$ contains a section,
which we call $\sigma_{\infty}$, and at most a
finite number of other points $z_i$.
\end{Prop}

\begin{Prop}\cite{ke00}
If $E \cong \PP(T_X^* \mid_x)$ is the exceptional divisor
of the blow-up $\tilde{X}  \to X$, then
the restricted morphism
$$\tilde{i}_x\mid_{\sigma_{\infty}}: \sigma_{\infty} \to E$$
is finite.
\end{Prop}

\begin{Prop}[=\cite{ke01} Proposition 3.1]
\hfill
\begin{enumerate}
\item
The map $$\tilde{i}_x\mid{\sigma_{\infty}}:\sigma_{\infty} \cong
\tilde{H}_x \to E$$ is an embedding. In particular, $\tilde{H}_x$ is
smooth.
\item The tangent map $T\tilde{i}_x$ has maximal rank
along $\sigma_{\infty}$.
 In particular, $N_{ \sigma_{\infty},\;U_{x} }
 \cong N_{E,\; \tilde{X} } \cong \mO_{\PP^{n-1}}(-1)$.
\end{enumerate}
\end{Prop}

We now sketch the proof that $X \cong \PP^n$. The argument  is taken
from \cite{ke01}. We do not claim any originality. Consider the
morphism $i_x: U_x \to X$ and its Stein factorization
$$\begin{CD}
 U_x  @>{\operatorname{\alpha}}>> Y @>{\operatorname{\beta}}>>X
\end{CD}$$
 where
$\alpha$ contracts the divisor $\sigma_{\infty}$,
and $\beta$ is a finite map.

Note that ${\bf R}^{1}\pi_{x} \mO_{U_x}=0$.
Pushing forward the exact sequence
\[ 0 \to \mO_{U_{x}} \to \mO_{U_{x}}(\sigma_{\infty})
\to \mO_{U_{x}}(\sigma_{\infty})\mid_{\sigma_{\infty}} \to 0\] gives
the sequence
\[ 0 \to \mO_{\PP^{n-1}} \to \mE \to  \mO_{\PP^{n-1}}(-1)  \to 0\]
where $\mE$ is a vector bundle of rank $2$ and $U_x \cong
\PP(\mE^{*})$.  Since $$Ext^{1}_{\PP^{n-1}}(\mO_{\PP^{n-1}}(-1),
\mO_{\PP^{n-1}})=0,$$ the sequence
\[ 0 \to \mO_{\PP^{n-1}} \to \mE \to  \mO_{\PP^{n-1}}(-1)  \to 0\]
splits, and the bundle $U_x \cong \PP(\mO_{\PP^{n-1}}(-1) \oplus
\mO_{\PP^{n-1}})$. This implies that there exists $\alpha': U_x \to
\PP^n$ which contracts $ \sigma_{\infty} \cong \PP^{n-1}$. Since
$\alpha$ and $\alpha'$ contract the same curve class, we have
$\alpha \cong \alpha'$. Therefore $Y \cong \PP^n$. To conclude the
proof,
 we need to show that
$$\begin{CD}
Y\cong \PP^n  @>{\beta}>> X
\end{CD}$$
is an isomorphism. Since $Y \cong \PP^n \to X$ is finite (and
surjective) and $X$ has only isolated singularities, we may assume
$X$ is smooth in comparing the $K_Y$ and $\beta^{*}K_X$. We obtain
$-K_Y= \beta^{*}(-K_X)+ R$ where $R$ is the (effective) ramification
divisor.

We need the following general fact: Let $C_1, C_2 \subset Y$ be any two curves. If $[C_1]=[ C_2]$  in the Mori cone
$\overline{NE}(Y)$, then $\beta_{*}[C_1]=\beta_{*}[C_2]$ in the Mori
cone $\overline{NE}(X)$ by the projection formula.

Let $d$ be the degree
of $\beta: Y \cong \PP^n \to X$.
Let $C'$ be the preimage of a general
curve $[C] \in H_x$. Consider the restriction  morphism $\beta \mid_{C'}: C' \to C$.
It is a $d$-to-$1$ cover.
Since $\rho(\PP^n)=1$, $[C']=k [l]$ where $k$
is a positive integer and $[l]$ is the class of a general line in
$\PP^n$.
Consider the image of the line $l$ under $\beta$. Since
$$k
\beta_{*}[l]=\beta_{*}[C']=d[C] \in \overline{NE}(X),$$ it follows
that the classes $\beta_{*}[l]$ and $[C]$ lie in  the same ray $R
\subset \overline{NE}(X)$. Note  that $k=1$  since $[C]$ is a curve
with the minimal $-K_X$-degree in the extremal ray $R$. It also
follows that
the $\beta_{*}[l]= \beta_{*}[C']=d[C] \in \overline{NE}(X)$.

Consider
$$n+1=l \cdot (-K_{\PP^n})= l \cdot [\beta^* (-K_X)+R]
=\beta_{*}l \cdot (-K_X)+ l \cdot R$$
$$ =dC \cdot (-K_X)+ l \cdot R
\geq dC \cdot (-K_X)=d(n+1).$$ This proves  that $d=1$ and
concludes
the proof.

\section{Threefolds}\label{threefoldsec}
We start with a simple proposition.
\begin{Prop}\label{crepant}
Let $X$ be a normal $\QQ$-Gorenstein projective variety of $dim
\;X=n \geq 3$. Assume that $C \cdot (-K_X) \geq n+1$ for every curve
$C \subset X$.
 Assume that there is a crepant partial resolution
$\pi: Y \to X$ such that $Y$ has only isolated LCIQ singularities.
Then $X \cong \PP^n$.
\end{Prop}
\begin{proof}
Since $\pi: Y \to X$ is crepant,  $K_Y= \pi^{*} K_X$. Let $R$ be an
extremal ray which is not contracted by $\pi:Y \to X$. Let $[C] \in
R$ be any curve. By assumption  $C \cdot (-K_Y) \geq n+1$. By
Theorem~\ref{main}, we have  $Y \cong \PP^n$. Note the Picard number
$\rho(\PP^n)=1$ and $X$ is normal. If $\pi:Y \cong \PP^n \to X$ is
not an isomorphism, then $\pi$ contracts at least one (and hence
all)  curve classes; that is, $X$ has to be a point. This is
impossible.
\end{proof}


 
\begin{proof}[Proof of Theorem~\ref{threefold}]
Take  any resolution $Y \to X$ and run the relative MMP over $X$.
Denote the resulting variety by $h: Z \to X$, where $Z$ is $\QQ$-factorial, has only terminal singularities (hence isolated LCIQ singularities), and $K_Z$ is $h$-nef. Write $K_Z=h^{*}K_X +\sum a_i E_i$ where $E_i$'s are exceptional. Since $K_Z$ is $h$-nef, it is clear that $a_i \leq 0$.

Case I: The morphism $h$ is crepant. 
In this case the Theorem follows from Proposition \ref{crepant}.

Case II: The morphism is not crepant, i.e. $K_Z=h^{*}K_X+ \sum_{i} a_i E_i$ with some $a_i<0$, say $a_1<0$.  Choose a curve $C\subset Z$ such that $C\cdot E_1>0$ and $C$ is not contained in any $E_i$. By the cone theorem for $Z$, we may write $[C]=\sum_j b_j v_j$ with $\mathbb{R}_{\geq 0}v_j$ an extremal ray in $\overline {NE}(Z)$ and $b_j>0$. Since $C\cdot \sum_i (-a_i E_i)\geq C\cdot (-a_1 E_1)>0$,  there is at least one $v_k$, say $v_1$, such that $v_1\cdot \sum_i (-a_i E_i) >0$. Hence 
\begin{equation}\label{nef}
v_1\cdot K_Z=v_1\cdot (h^*K_X+\sum_i a_i E_i)< 0,
\end{equation}
and therefore there is a curve class which generates the extremal ray $\RR_{\geq 0}v_1$. Let $C_1$ be any curve whose class is in the ray $\RR_{\geq 0}v_1$. We claim that $C_1$ is not contracted by $h$: suppose otherwise, then $C_1\cdot K_Z\geq 0$ since $K_Z$ is $h$-nef; this contradicts (\ref{nef}). 

Since $h_*C_1\cdot (-K_X)\geq 4$ by assumption and $C_1\cdot \sum_i (-a_i E_i)>0$, it follows that $$C_1\cdot (-K_Z)=C_1\cdot h^*(-K_X)+C_1\cdot \sum_i (-a_i E_i)$$ $$=h_*C_1\cdot (-K_X)+C_1\cdot \sum_i (-a_i E_i)\geq 4.$$ By Theorem \ref{main}, $Z\simeq \PP^3$. 

Since $\text{Pic}(Z)= \text{Pic}(\PP^3)= \ZZ$ and $h:Z\simeq \PP^3\to X$ is birational, we conclude that $h$ is an isomorphism and $X\simeq \PP^3$.
\end{proof}

\end{document}